\documentclass[10pt]{amsart}
\usepackage{geometry}                
\usepackage{graphicx}
\usepackage{dsfont}
\usepackage{amssymb}
\usepackage{amsmath}
\usepackage{epstopdf}
\usepackage{fullpage}
\usepackage{enumerate}
\usepackage{color,subcaption}
\usepackage{amsthm}
\usepackage{placeins}
\usepackage{gensymb} 
\usepackage{array}
\usepackage{blkarray}
\usepackage{multirow}
\usepackage{float}
\usepackage{morefloats}
\usepackage{pgfplots}
\usepackage[abs]{overpic}
\usepackage{tikz}
\usetikzlibrary{calc,positioning}

\usepackage{amsopn}

\usepackage[colorlinks=true]{hyperref}
\usepackage[nameinlink,capitalise]{cleveref}
\crefname{equation}{}{}

\makeatletter
\def\subsection{\@startsection{subsection}{3}%
  \z@{.5\linespacing\@plus.7\linespacing}{.5\linespacing}%
  {\bf}}
\makeatother

\makeatletter
\def\subsubsection{\@startsection{subsubsection}{3}%
  \z@{.5\linespacing\@plus.7\linespacing}{.5\linespacing}%
  {\it}}
\makeatother

\let\oldtocsection=\tocsection

\let\oldtocsubsection=\tocsubsection

\let\oldtocsubsubsection=\tocsubsubsection

\renewcommand{\tocsection}[2]{\hspace{0em}\oldtocsection{#1}{#2}\textbf}
\renewcommand{\tocsubsection}[2]{\hspace{1em}\oldtocsubsection{#1}{#2}}
\renewcommand{\tocsubsubsection}[2]{\hspace{2em}\oldtocsubsubsection{#1}{#2}}

\newcommand{\dv}{\text{\rm div}}

\renewcommand{\d}{\text{\rm d}}

\newcommand{\e}{\varepsilon}

\newcommand{\calP}{{\mathcal P}}

\newcommand{\calC}{{\mathcal C}}

\newcommand{\calD}{{\mathcal D}}

\newcommand{\Vol}{\text{\rm Vol}}

\newcommand{\R}{{\mathbb R}}
\renewcommand{\P}{{\mathbb P}}
\newcommand{\Q}{{\mathbb Q}}

\newcommand{\Uad}{{\mathcal U}_{\text{\rm ad}}}
\newcommand{\CD}[1]{{#1}}

\usepackage{soul,xcolor}
\setstcolor{violet}

\begin{document}
\newtheorem{theorem}{Theorem}[section]
\newtheorem{problem}{Problem}[section]
\newtheorem{remark}{Remark}[section]
\newtheorem{example}{Example}[section]
\newtheorem{definition}{Definition}[section]
\newtheorem{lemma}{Lemma}[section]
\newtheorem{corollary}{Corollary}[section]
\newtheorem{proposition}{Proposition}[section]
\numberwithin{equation}{section}

\title{Entropy-regularized Wasserstein distributionally robust shape and topology optimization}
\author{
C. Dapogny$^1$, F. Iutzeler$^1$, A. Meda$^1$ and B. Thibert $^1$
}

\maketitle
\begin{center}
\emph{\textsuperscript{1} Univ. Grenoble Alpes, CNRS, Grenoble INP\footnote{Institute of Engineering Univ. Grenoble Alpes}, LJK, 38000 Grenoble, France.}
\end{center}
 
\begin{abstract}
This brief note aims to introduce the recent paradigm of distributional robustness in the field of shape and topology optimization. 
Acknowledging that the probability law of uncertain physical data is rarely known beyond a rough approximation constructed from observed samples,
we optimize the worst-case value of the expected cost of a design when the probability law of the uncertainty is ``close'' to the estimated one up to a prescribed threshold. 
The ``proximity'' between probability laws is quantified by the Wasserstein distance, \CD{a notion pertaining} to optimal transport theory. \CD{The combination of the} classical entropic regularization technique in this field
with recent results from convex duality theory allows to reformulate the distributionally robust optimization problem in a way which is tractable for computations.
Two numerical examples are presented, in the different settings of density-based topology optimization and geometric shape optimization. 
They exemplify the relevance and applicability of the proposed formulation regardless of the \CD{selected} optimal design framework.
\end{abstract} 

\bigskip
\bigskip
\hrule
\tableofcontents
\vspace{-0.5cm}
\hrule
\bigskip
\bigskip

\section{Introduction}\label{sec.intro}

\noindent In realistic situations, the physical behavior of a design $h$ depends on one or several parameters, collectively denoted by $\xi\in \Xi$; 
for instance, when $h$ is a mechanical structure, $\xi$ may account for the loads applied on $h$ or the coefficients of the constituent elastic material. 
These parameters $\xi$ affect, often dramatically, the physical response $u_{h,\xi}$ and the cost $\calC(h,\xi)$ of $h$, which raises the need to incorporate a degree of awareness to uncertainties over their values 
in optimal design procedures \cite{cherkaev1999optimal,maute2014topology}. 

The treatment of uncertainties in optimal design usually fits in one of the following two frameworks. 
When no information is available about the uncertain parameters $\xi$ except for a maximum bound $m$ on their amplitude $|| \xi ||$,
\CD{the only realistic resort is a} worst-case design approach, \CD{which consists} in minimizing the \CD{largest value} $h \mapsto \sup_{|| \xi ||\leq m}\calC(h,\xi)$ \CD{of $\calC(h,\xi)$ among all possible scenarii}. 
When the law $\P_{\text{true}}$ of $\xi$ is known, probabilistic approaches are usually preferred: 
the expectation $h \mapsto \int_\Xi \calC(h,\xi) \: \d \P_{\text{true}}(\xi)$ (or another risk measure) of the cost can be minimized, see e.g. \cite{shapiro2021lectures}.
Both approaches suffer from major drawbacks: while worst-case \CD{settings} are often deemed too ``pessimistic'', 
insofar as the robust design may show poor nominal performance for the sake of providing for an unlikely worst-case scenario, 
the law $\P_{\text{true}}$ of $\xi$ featured in probabilistic approaches is often unknown, and in practice chosen in a heuristic manner.

Recently, the paradigm of distributionally robust optimization has emerged as an elegant means to overcome this conceptual shortcoming of probabilistic formulations, 
see \cite{lin2022distributionally,rahimian2019distributionally,zhen2021mathematical}.
Building on a nominal probability law $\P$ for $\xi$, which is for instance constructed from observed data,
it minimizes the worst value $\sup_{\Q} \int_\Xi \calC(h,\xi) \:\d \Q(\xi)$ of the expected cost when the law $\Q$ of $\xi$ is ``close'' to $\P$ within a given tolerance $m$.

Hitherto, this idea has been considered mainly in academic settings; the purpose of this note is to show how it can be applied in the field of optimal design.
More precisely, we adopt the viewpoint in \cite{gao2016distributionally,mohajerin2018data} where the notion of ``closeness'' between \CD{probability} measures is quantified by the Wasserstein distance from optimal transport theory \cite{merigot2021optimal,santambrogio2015optimal}. We then take advantage of the key results from convex duality proved in \cite{azizian2022regularization,wang2021sinkhorn} 
to reformulate the entropy-regularized version of the distributionally robust optimal design problem in a manner which is amenable to computations.

In principle, this methodology can be implemented in any optimal design framework. 
For simplicity, the main ideas are presented in a formal and non technical way in the context of a model density-based topology optimization problem, 
but we also propose a numerical example in the setting of (geometric) shape optimization.
A longer article, containing full mathematical details and extensive numerical experiments, is currently in preparation about this topic. 

\section{Presentation of the distributionally robust optimal design problem}\label{sec.desc}

\subsection{The deterministic compliance minimization problem}\label{sec.detpb}

\noindent The considered designs are elastic structures contained in a fixed hold-all domain $D \subset \R^d$, 
that are clamped on a region $\Gamma_D \subset \partial D$ and subjected to traction loads on a disjoint subset $\Gamma_N \subset \partial D$.
They are represented as density functions $h$ on $D$, i.e. 
$$ h \in \Uad, \quad \Uad := L^\infty(D,[0,1]),$$
\CD{with the meaning that} $h(x)$ equals $0$ (resp. $1$) at points $x \in D$ surrounded by void (resp. by material) 
and \CD{ that $h(x)$ takes values in $(0,1)$ in ``grayscale'' regions of $D$}, made of a fictitious mixture of material and void. 
The material properties \CD{induced in $D$ by such a density function} are encoded in the Hooke's tensor $A(h)$, 
which is related to $h$ via the so-called SIMP law:
\begin{equation}\label{eq.simplaw}
A(h)(x) = \Big(\eta + (1-\eta)h(x)^p\Big)A, \quad x \in D, 
\end{equation}
where $A$ is the Hooke's law of the reference elastic material, and $\eta \ll 1$ is a small parameter mimicking the presence of void, 
see e.g. \cite{bendsoe2013topology} about this classical setting.

The uncertain parameters \CD{at play} are the loads applied on $\Gamma_N$.
These are assumed to be \CD{given by a constant vector $\xi$} which belongs to a ``large enough'' closed ball $\Xi = \overline{B(0,R)} \subset \R^d$, for some $R >0$. 
The displacement of a design $h \in \Uad$ in response to the load $\xi \in \Xi$ is then the solution $u_{h,\xi}:D \to \R^d$ to the linear elasticity system:
\begin{equation}
\left\{
\begin{array}{cl}
-\dv (A(h) e(u_{h,\xi})) = 0 & \text{in } D,\\
u_{h,\xi} = 0 & \text{on } \Gamma_D, \\ 
A(h)e(u_{h,\xi}) n=\xi & \text{on } \Gamma_N, \\ 
A(h)e(u_{h,\xi}) n=0 & \text{on } \partial D \setminus (\overline{\Gamma_D} \cup \overline{\Gamma_N}),
\end{array}
\right.
\end{equation}
where $e(u) := \frac12 (\nabla u + \nabla u^T)$ denotes the strain tensor associated to a displacement field $u : D \to \R^d$. 

In this \CD{situation}, when the load $\xi \in \Xi$ is known exactly, the optimization problem of interest reads
\begin{equation} \label{eq.nomoptpb}
\min\limits_{h \in \Uad} \calC(h,\xi) \:\:  \text{ s.t. } \:\: \Vol(h) = V_T.
\end{equation}
Here, the cost $\calC(h,\xi)$ of a design $h \in \Uad$ submitted to the load $\xi \in \Xi$ is the compliance, i.e.
$$ \calC(h,\xi) = \int_D A(h)e(u_{h,\xi} ) : e(u_{h,\xi}) \:\d x,$$
the volume functional is \CD{defined} by  $\Vol(h) = \int_D h \:\d x$, and $V_T$ is a volume target.

\subsection{The distributionally robust compliance minimization problem}

\noindent We now turn to the more realistic situation where the load $\xi$ applied on $\Gamma_N$ is governed by a probability law $\Q$ on $\Xi$ which is unknown. 
Fortunately, in most practical situations, $\Q$ can be estimated by a nominal law $\P$, which is typically the empirical sum of a collection of observations $\xi_i \in \Xi$, $i=1,\ldots,N$:
$$ \P = \frac1N \sum\limits_{i=1}^N \delta_{\xi_i}.$$

Let $\calP(\Xi)$ be the space of probability measures on the compact set $\Xi \subset \R^d$.
\CD{We equip the latter} with the Wasserstein distance stemming from optimal transport theory, \CD{which is} defined by: 
\begin{equation}\label{eq.Wass}
 W(\P,\Q) = \inf \int_{\Xi \times \Xi} c(\xi,\zeta) \: \d \pi(\xi,\zeta).
 \end{equation}
In this formula, \CD{the ``ground cost'' $c : \Xi \times \Xi \to \R$ is a continuous function between loads in $\Xi$, 
which is throughout the sequel chosen as the squared Euclidean distance $c(\xi,\zeta) := |\xi-\zeta|^2$}.
The infimum in \cref{eq.Wass} is taken over all transport plans \CD{-- or ``couplings'' --} $\pi \in \calP(\Xi \times \Xi)$ between $\P$ and $\Q$. 
 These are probability measures on the product space $\Xi \times \Xi$ 
whose first and second marginals $\pi_1$ and $\pi_2$ coincide with $\P$ and $\Q$, respectively, \CD{that is, for all continuous functions $f$ and $g$ on $\Xi$, it holds:
$$ \int_{\Xi \times \Xi} f(\xi) \:\d\pi(\xi,\zeta)  = \int_{\Xi} f(\xi) \:\d \P(\xi), \text{ and }   \int_{\Xi \times \Xi} g(\zeta) \:\d\pi(\xi,\zeta)  = \int_{\Xi} g(\zeta) \:\d \Q(\zeta).$$
}
Intuitively, in \cref{eq.Wass}, $c(\xi,\zeta) $ measures the cost of ``moving a unit of mass'' from $\xi$ to $\zeta$, 
 $\pi(\xi,\zeta)$ encodes the ``quantity of mass'' transported from $\xi$ to $\zeta$, \CD{so that $W(\P,\Q)$ evaluates the minimum cost of moving the mass from $\P$ onto that of $\Q$ among all transport plans}. 
We refer e.g. to \cite{santambrogio2015optimal} about the properties of this distance,
and to \cite{peyre2019computational} for an overview of its use in applications. \par\medskip

With these definitions at hand, the distributionally robust counterpart of \cref{eq.nomoptpb} is
\begin{equation}\label{eq.dropb}
 \min \limits_{h \in \Uad} J_{\text{dr}}(h) \:\: \text{ s.t. } \:\: \Vol(h) = V_T,
 \end{equation}
where $J_{\text{dr}}(h)$ is the worst (maximum) value of the expected cost when the law $\Q \in \calP(\Xi)$ of the uncertain parameter $\xi$ 
is at distance less than a given threshold $m$ from the nominal law $\P$:
\begin{equation}\label{eq.Jdr}
 J_{\text{dr}}(h) = \sup\limits_{\Q \in \calP(\Xi) \atop W(\P,\Q) \leq m} \int_\Xi  \calC(h,\xi) \:\d\Q(\xi).
 \end{equation}

\section{Entropic regularization of the distributionally robust problem}\label{sec.ref}

\noindent The distributionally robust optimal design problem \cref{eq.dropb}
is hard to tackle as is. To alleviate this issue,
we consider the entropy-regularized version of the Wasserstein distance proposed in \cite{cuturi2013sinkhorn}:
\begin{equation}\label{eq.We}
W_\e(\P,\Q) = \inf\limits_{\pi \in \calP(\Xi \times \Xi) \atop \pi_1 = \P, \: \pi_2 = \Q} \left\{ \int_{\Xi \times \Xi} c(\xi,\zeta) \:\d \pi(\xi,\zeta) + \e H(\pi) \right\}, 
\end{equation}
where $\e>0$ is a ``small'' smoothing parameter.
\CD{The definition of the entropy $H(\pi)$ of an element $\pi \in \calP(\Xi \times \Xi)$ involves a fixed reference coupling $\pi_0 \in \calP(\Xi \times \Xi)$ playing the role of a ``prior''; it accounts for the deviation of $\pi$ from $\pi_0$ as:}
$$ H(\pi) = \left\{
\begin{array}{cl}
\int_{\Xi \times \Xi} \log \frac{\d \pi}{\d \pi_0} \:\d \pi & \text{if } \pi \text{ is absolutely continuous w.r.t. } \pi_0, \\
\infty & \text{otherwise.}
\end{array}
\right. $$
According to \cite{azizian2022regularization}, a judicious choice about $\pi_0$ is provided by the following formula:
\begin{equation}\label{eq.pi0}
 \pi_0(\xi,\zeta) =  \P(\xi) \d \nu_\xi(\zeta),  \:\: \text{ with } \:\: \d\nu_\xi(\zeta) :=  \alpha_\xi e^{-\frac{c(\xi,\zeta)}{2\sigma}}  \mathds{1}_{\Xi}(\zeta)  \d \zeta ,
 \end{equation}
for some $\sigma >0$ and a normalization factor $\alpha_\xi$ ensuring that $\d\nu_\xi$ is a probability distribution on $\Xi$.
Precisely, the above definition means that
$$ \text{For all continuous functions } \varphi :\Xi \times \Xi \to \R,  \quad \int_{\Xi \times \Xi} \varphi(\xi,\zeta)\: \d\pi_0(\xi,\zeta) = \int_\Xi \left( \int_\Xi \varphi(\xi,\zeta) \: \d \nu_\xi(\zeta) \right) \:\d \P(\xi).$$

\CD{
Without entering into precise statements, let us mention that this choice of $\pi_0$ guarantees explicit bounds in terms of $\e$ between the non regularized and regularized quantities $\sup_{W(\P,\Q) \leq m} \int_\Xi f(\xi) \:\d \Q(\xi)$ and $\sup_{W_\e(\P,\Q) \leq m} \int_\Xi f(\xi) \:\d \Q(\xi)$ attached to an arbitrary continuous function $f$ on $\Xi$, see \cite{azizian2022regularization}.
}\par\medskip

\CD{With these notions at hand, we} introduce the entropy-regularized version of the problem \cref{eq.dropb}:
\begin{equation}\label{eq.dropbeps}
 \min \limits_{h \in \Uad} J_{\text{dr},\e}(h) \:\: \text{ s.t. } \:\: \Vol(h) = V_T, \text{ where } J_{\text{dr},\e}(h) := \sup\limits_{\Q \in \calP(\Xi) \atop W_\e(\P,\Q) \leq m} \int_\Xi  \calC(h,\xi) \:\d\Q(\xi).
 \end{equation}
This program is intricate at first glance, as it features nested maximization and minimization problems. 
Fortunately, the functional $J_{\text{dr},\e}(h) $ admits a convenient dual reformulation as \CD{an infimum}, as expressed by the following result from convex analysis. 
\CD{The latter is proved rigorously in a more general context in \cite{azizian2022regularization} (see also \cite{wang2021sinkhorn});
 for the convenience of the reader, we provide an intuitive sketch of the main arguments in \cref{app.proof}.} 

\begin{proposition}\label{prop.supWass}
Let $\Xi$ be a convex and compact subset of $\R^d$, $f : \Xi \to \R$ be a continuous function and let $\P \in \calP(\Xi)$ be a probability measure. 
\CD{For any $m >0$,} and for a sufficiently small value of $\sigma$, the following equality holds:
\begin{equation}\label{eq.dualform}
\sup\limits_{W_\e(\P,\Q) \leq m}  \int_\Xi f(\xi) \:\d \Q(\xi) = 
\inf \limits_{\lambda \geq 0 } \left\{ \lambda m + \lambda \e \int_\Xi \log \left( \int_\Xi e^{\frac{f(\zeta) - \lambda c(\xi,\zeta) }{\lambda \e} } \d \nu_\xi(\zeta) \right) \d \P(\xi) \right\}.
\end{equation}
\end{proposition}\par\medskip

Taking advantage of this result, the distributionally robust optimization problem \cref{eq.dropbeps} rewrites:
\begin{equation}\label{eq.droaug}
 \min \limits_{h \in \Uad \atop \lambda \geq 0} \calD(h,\lambda) \text{ s.t. } \Vol(h) = V_T,
 \end{equation}
\begin{equation}\label{eq.calDhl}
\text{where } \calD(h,\lambda) := \lambda m 
 +  \lambda \e \int_\Xi \log \left( \int_\Xi e^{\frac{\calC(h,\zeta) - \lambda c(\xi,\zeta)}{\lambda \e}} \d \nu_\xi(\zeta) \right) \d \P(\xi).
 \end{equation}
 This new version \cref{eq.droaug} boils down to a single minimization problem for the pair $(h,\lambda)$. It can be solved by a standard constrained optimization algorithm
 \CD{based on} the derivatives of the objective functional $\calD(h,\lambda)$ with respect to both variables $h$ and $\lambda$. \CD{The calculation of the latter}
follows from a standard (albeit a little tedious) adjoint-based procedure, \CD{see e.g. \cite{allaire2020survey,allaire2007conception,bendsoe2013topology}.}

\begin{remark}
A duality result similar to \cref{prop.supWass} actually holds when the regularized quantity $W_\e(\P,\Q)$ is replaced by the true Wasserstein distance $W(\P,\Q)$,
thus leading to a reformulation of the distributionally robust problem \cref{eq.dropb} of the form \cref{eq.droaug}. The latter is however more difficult to handle from the numerical viewpoint
since it involves the supremum $ \sup_{\zeta \in \Xi}\Big( \calC(h,\zeta) - \lambda c(\xi,\zeta)\Big)$ in place of the ``smooth'' log-sum-exp approximation $ \e \log \left( \int_\Xi e^{\frac{\calC(h,\zeta) - \lambda c(\xi,\zeta)}{\lambda \e}} \d \nu_\xi(\zeta) \right) $ featured in \cref{eq.calDhl}, see \cite{gao2016distributionally}.
\end{remark}

\section{Numerical results}\label{sec.num}

\subsection{Topology optimization of a 2d bridge}\label{sec.br}

\begin{figure}[!ht]
\centering
\begin{tabular}{ccc}
\begin{minipage}{0.275\textwidth}
\vspace{-0.3cm}
\hspace{-0.3cm}
\includegraphics[width=1.0\textwidth]{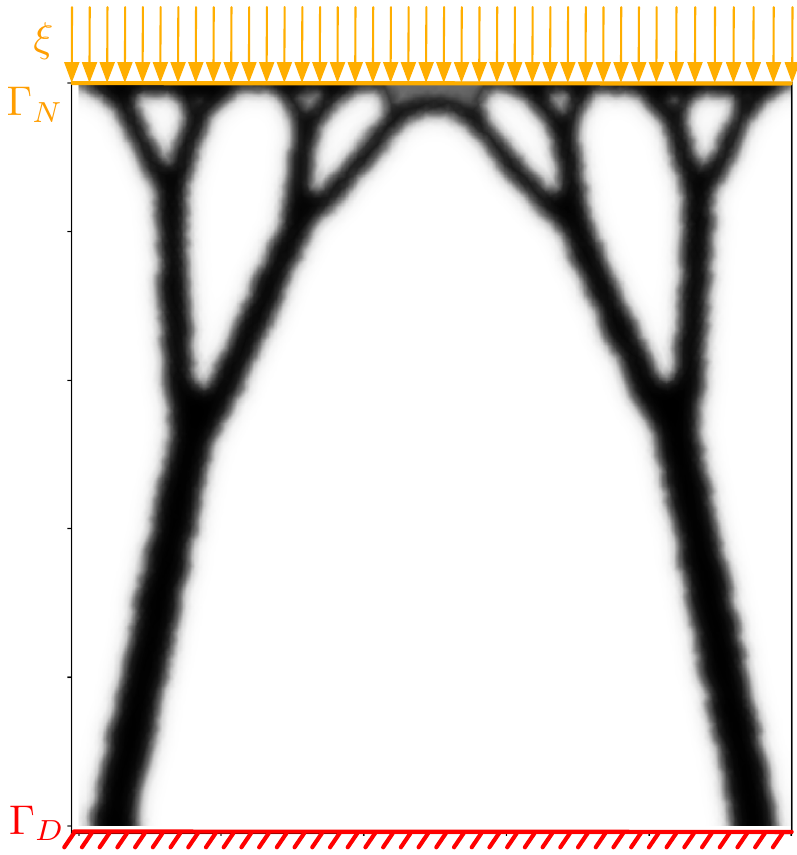}
\end{minipage} &
\begin{minipage}{0.25\textwidth}
\includegraphics[width=1.0\textwidth]{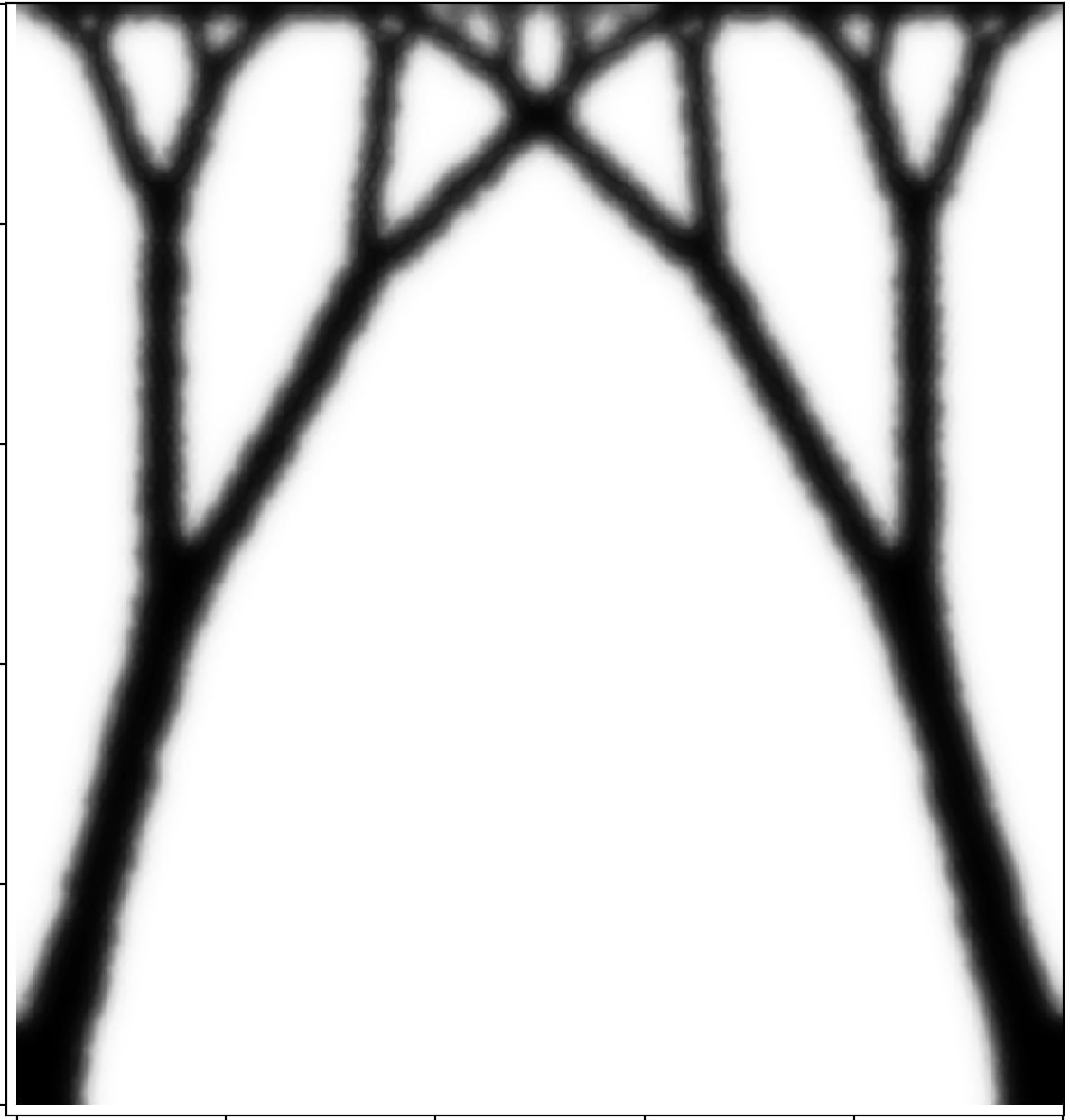}
\end{minipage}&
\begin{minipage}{0.25 \textwidth}
\includegraphics[width=1.0\textwidth]{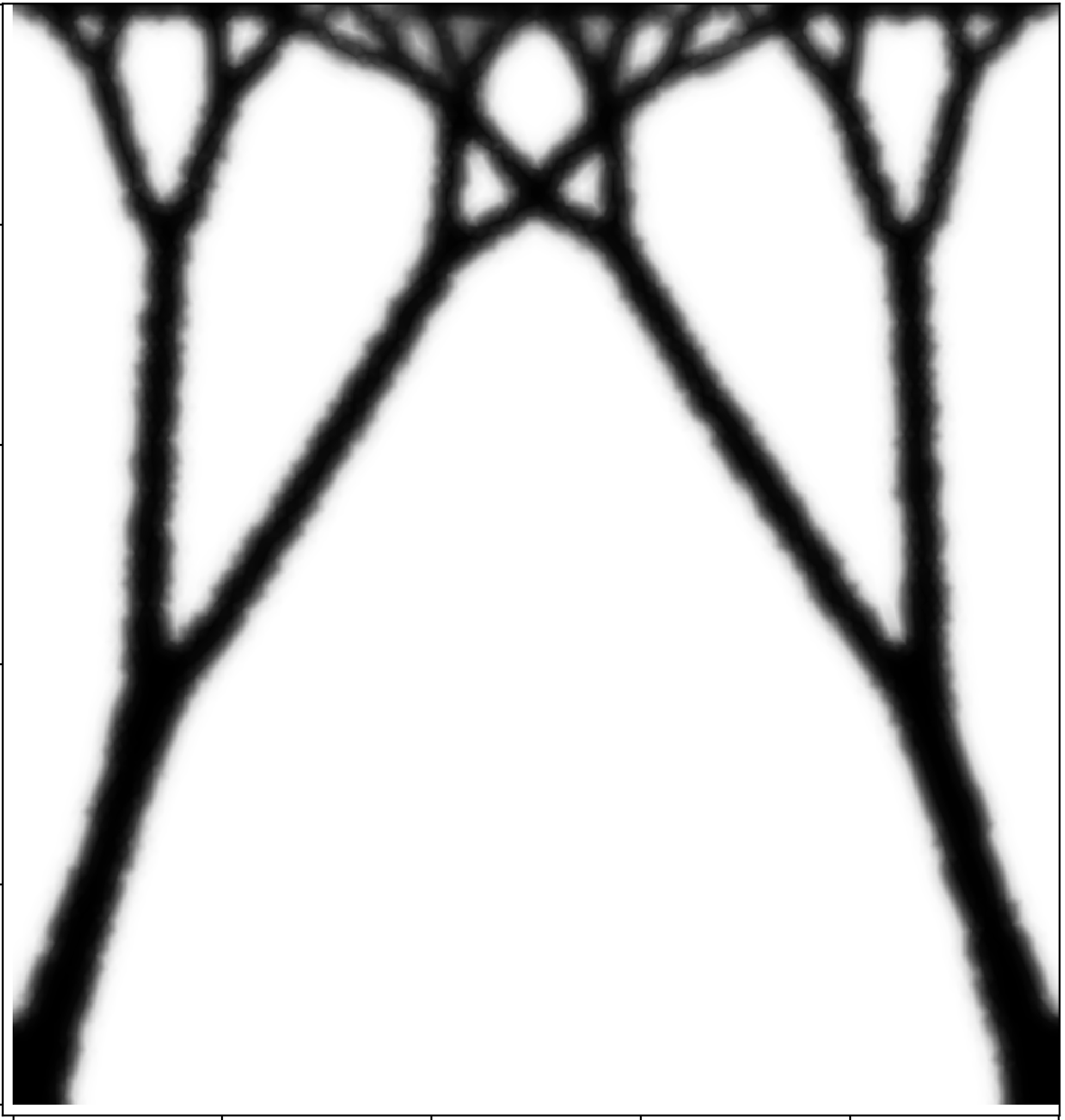}
\end{minipage} \\ 
\begin{minipage}{0.25\textwidth}
\includegraphics[width=1.0\textwidth]{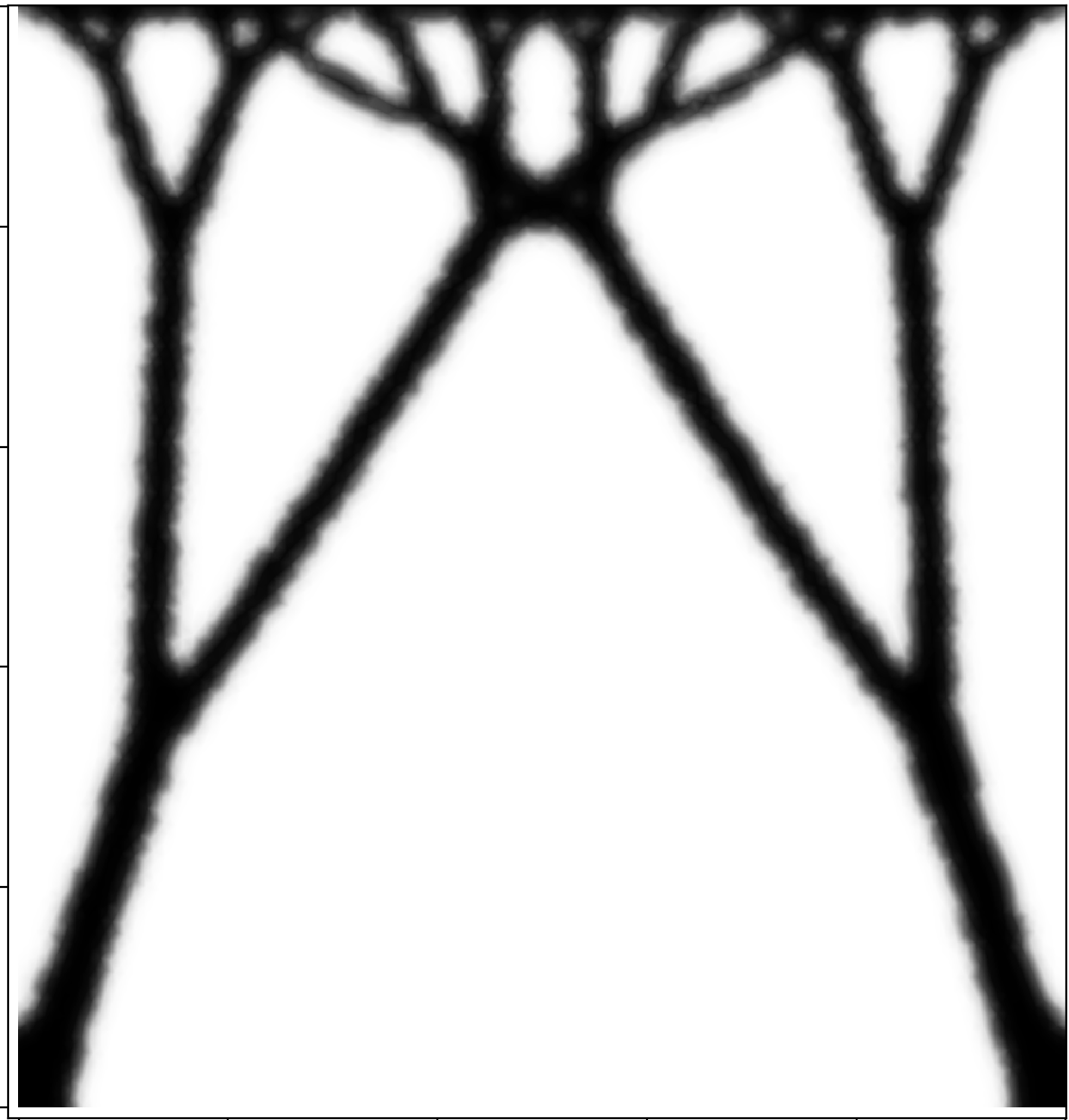}
\end{minipage}  &
\begin{minipage}{0.25\textwidth}
\includegraphics[width=1.0\textwidth]{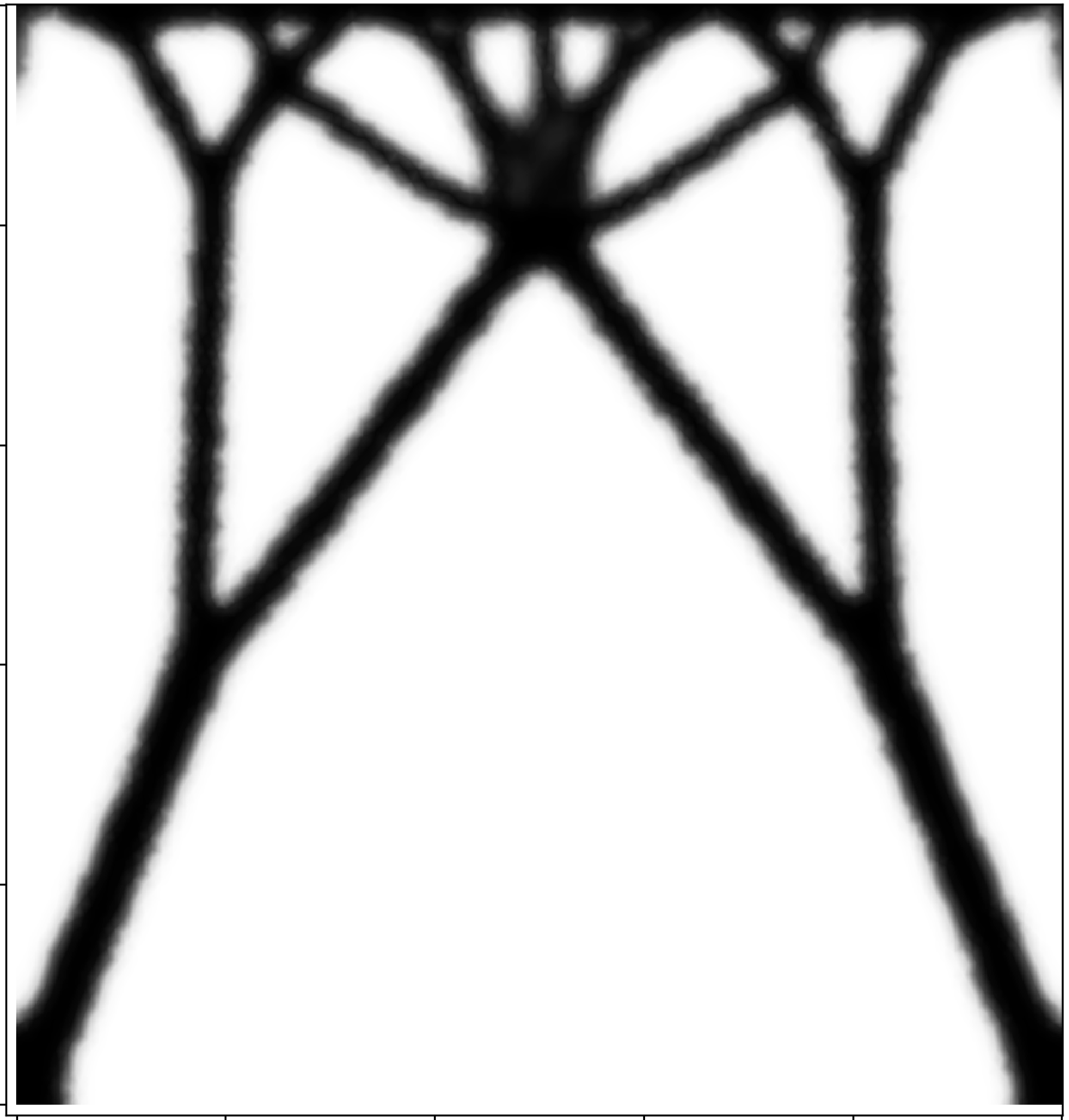}
\end{minipage} &
\begin{minipage}{0.25\textwidth}
\includegraphics[width=1.0\textwidth]{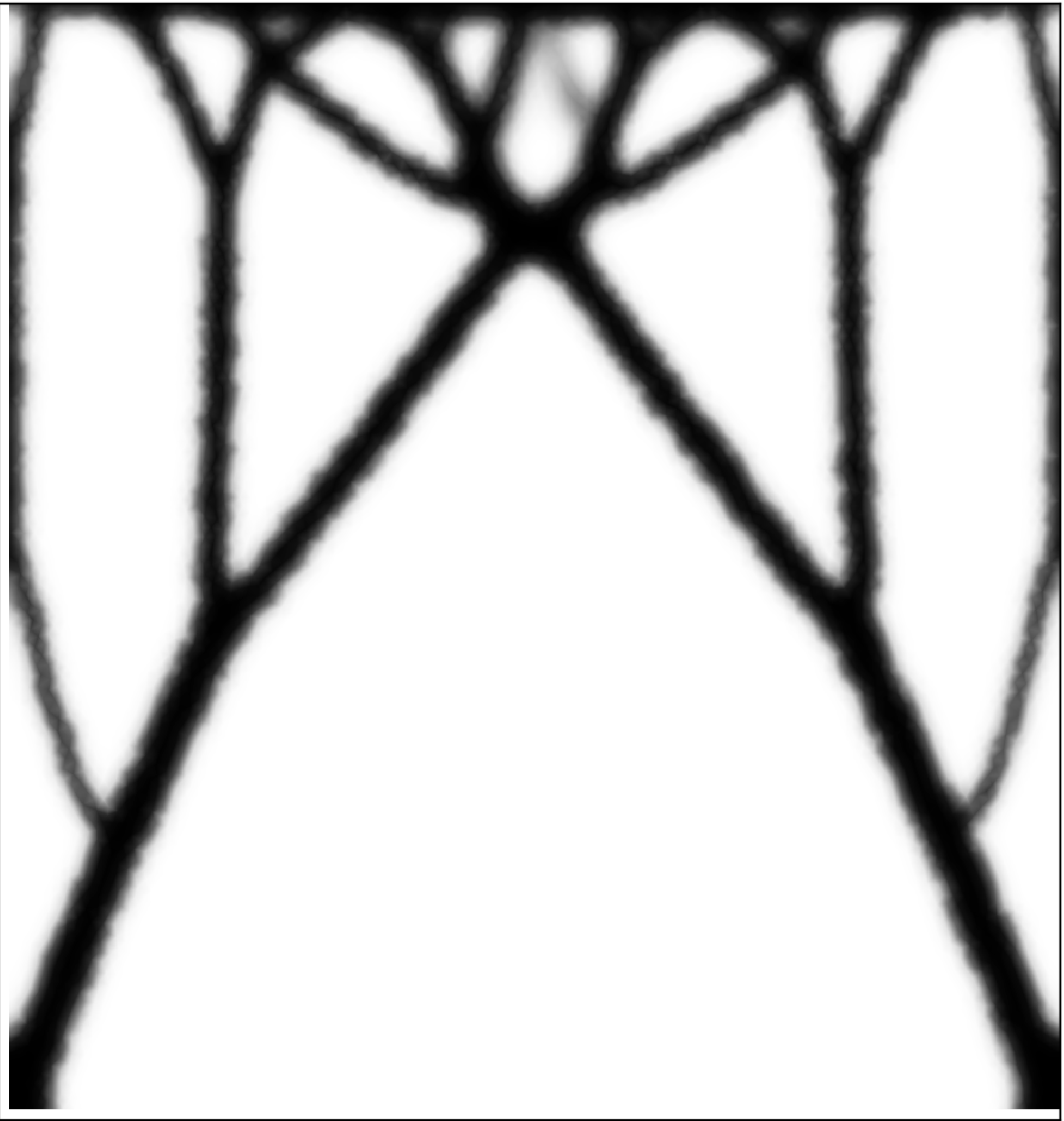}
\end{minipage} 
\end{tabular}

\caption{\it (From left to right, top to bottom) Optimized density in the bridge topology optimization example of \cref{sec.br} for $m=0$ (with details of the test-case) and $m=0.25$, $0.52$, $0.6$, $0.9$, $1$.}
\label{fig.br}
\end{figure}

\noindent Our first numerical example unfolds in the density-based context of \cref{sec.detpb}, 
and deals with the topology optimization of a 2d bridge. The considered designs $h$ are contained in a box $D$ with size $1 \times 1$; 
they are clamped on the bottom side $\Gamma_D$ of $\partial D$ and subjected to a constant load $\xi$ applied on the whole upper side $\Gamma_N$, see \cref{fig.br} (top, left) for details.
The nominal probability law $\P$ for $\xi$ is constructed from one single observation $\xi_1$, corresponding to a unit vertical load:
$$ \P = \delta_{\xi_1}, \quad \xi_1 = (0,-1).$$
The entropic regularization coefficient $\e$ is $1e^{-2}$ 
and the parameter $\sigma$ appearing in the reference coupling $\pi_0$ in \cref{eq.pi0} equals $1e^{-3}$.

\begin{figure}[!ht]
\centering
\includegraphics[width=0.8\textwidth]{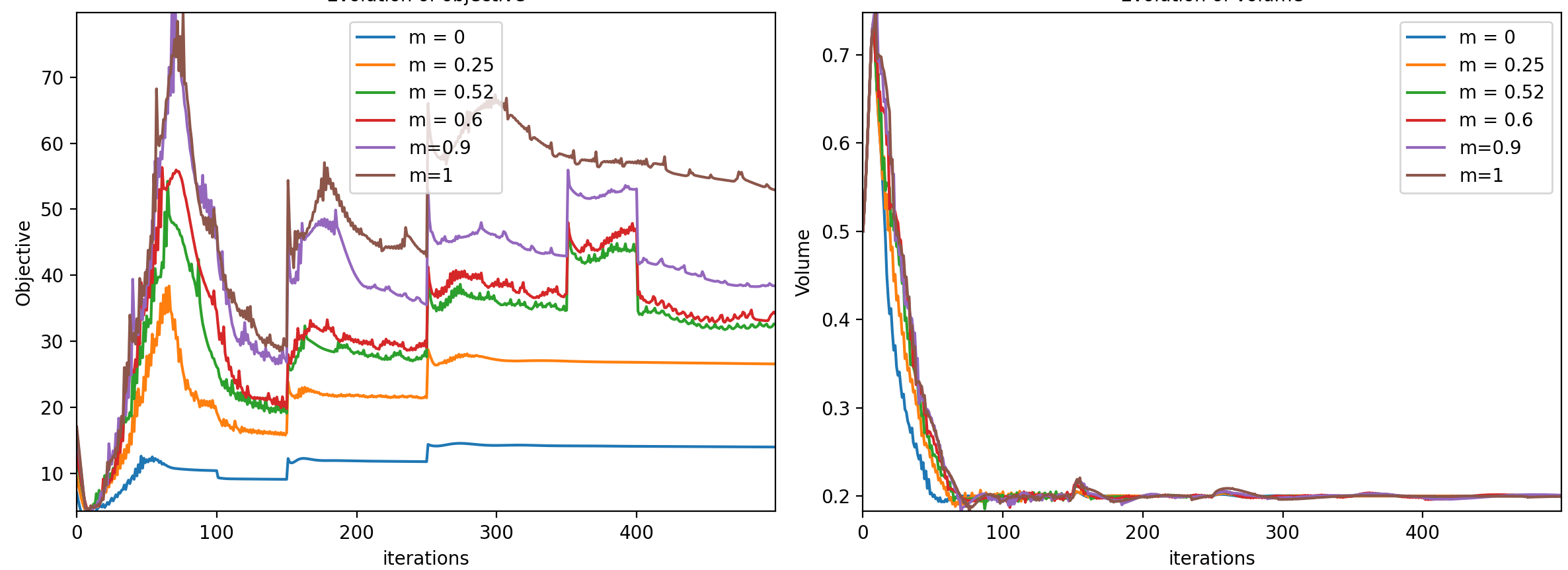}

\caption{\it Convergence histories for the experiments conducted with the bridge topology optimization example of \cref{sec.br}; the large jumps in the values of the objective function at particular iterations are due to an increase in the parameter $p$ of the material law \cref{eq.simplaw}.}
\label{fig.brhisto}
\end{figure}

The distributionally robust topology optimization problem \cref{eq.droaug} is solved with the target volume $V_T = 0.2$, for several values of the tolerance parameter $m$. 
The optimized designs are represented on \cref{fig.br}, and the histories of the computation are displayed on \cref{fig.brhisto}.
Understandably enough, the optimized designs develop thin branches to cope with larger loads with horizontal components,  
and their nominal performance $\calC(h,\xi_1)$ gets increasingly bad as $m$ grows, see \cref{tab.brtab}.

\begin{table}[!h]
\centering
	\begin{tabular}{|l|c|c|c|c|c|c|}
	\hline
	 Value of $m$  & $0$    &   $0.25$   & $0.52$   & $0.6$ & $0.9$ & $1$   \\
	\hline
	Nominal compliance & $13.9902$    & $17.3063$      & $19.2063$      &   $19.6829$   & $24.3765$ & $30.2474$       \\
	\hline
	\end{tabular}
	\caption{\it Values of the nominal cost for the optimized bridges of \cref{sec.br}.}
\label{tab.brtab}
\end{table}

\subsection{Geometric optimization of a 2d cantilever}\label{sec.canti}

\begin{figure}[!ht]
\centering
\begin{tabular}{cc}
\begin{minipage}{0.54\textwidth}
\vspace{0.1cm}
\hspace{-0.03cm}
\includegraphics[width=1.0\textwidth]{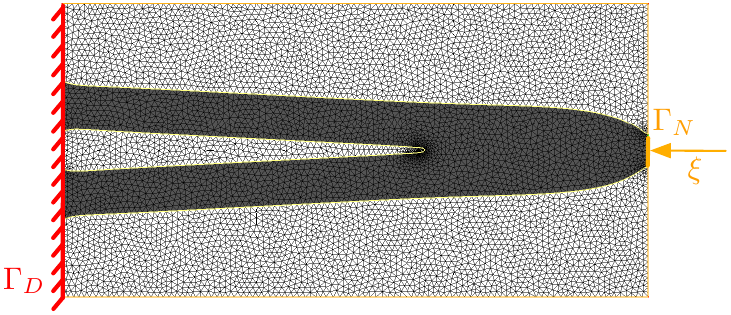}
\end{minipage} &
\begin{minipage}{0.43\textwidth}
\includegraphics[width=1.0\textwidth]{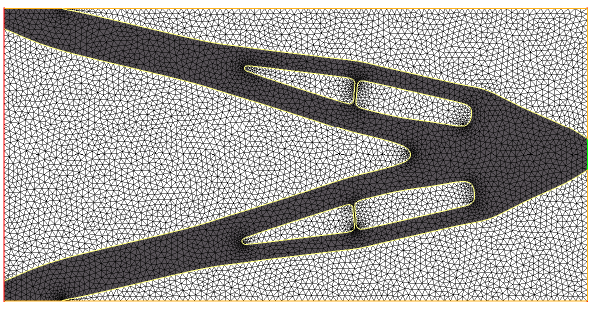}
\end{minipage} \\
\begin{minipage}{0.43\textwidth}
\includegraphics[width=1.0\textwidth]{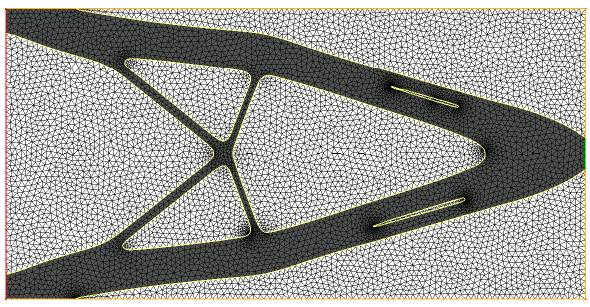}
\end{minipage} &
\begin{minipage}{0.43\textwidth}
\includegraphics[width=1.0\textwidth]{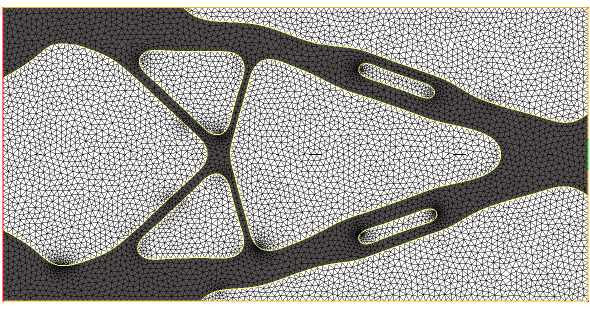}
\end{minipage} 
\end{tabular}

\caption{\it (From left to right, top to bottom) Optimized shape in the cantilever shape optimization example of \cref{sec.canti} for $m=0$ (with details of the test-case), optimized shapes for $m=1$, $1.5$, $2$.}
\label{fig.canti}
\end{figure}

\noindent Our second example is about the optimization of a 2d cantilever beam, 
and it is treated from the geometric shape optimization viewpoint \cite{allaire2007conception,allaire2020survey}: the considered designs are domains $\Omega$
contained in the fixed computational domain $D = [0,2] \times [0,1]$.
They are clamped on their left-hand side $\Gamma_D \subset \partial D$ 
and a constant load $\xi \in \Xi$ is applied on a small region $\Gamma_N$ at the middle of their right-hand side, see \cref{fig.canti} for the details.  
We rely on the mesh evolution method from our previous work \cite{allaire2014shape,dapogny2022tuto} to track the evolution of the mesh of the optimized shape. 

Again, the nominal law $\P$ for the load is constructed from only one sample $\xi_1$:
$$ \P = \delta_{\xi_1}, \quad \xi_1=(-1,0),$$
and the parameters $\e$ and $\sigma$ equal respectively $1e^{-2}$ and $1e^{-3}$.

\begin{figure}[!ht]
\centering
\includegraphics[width=0.8\textwidth]{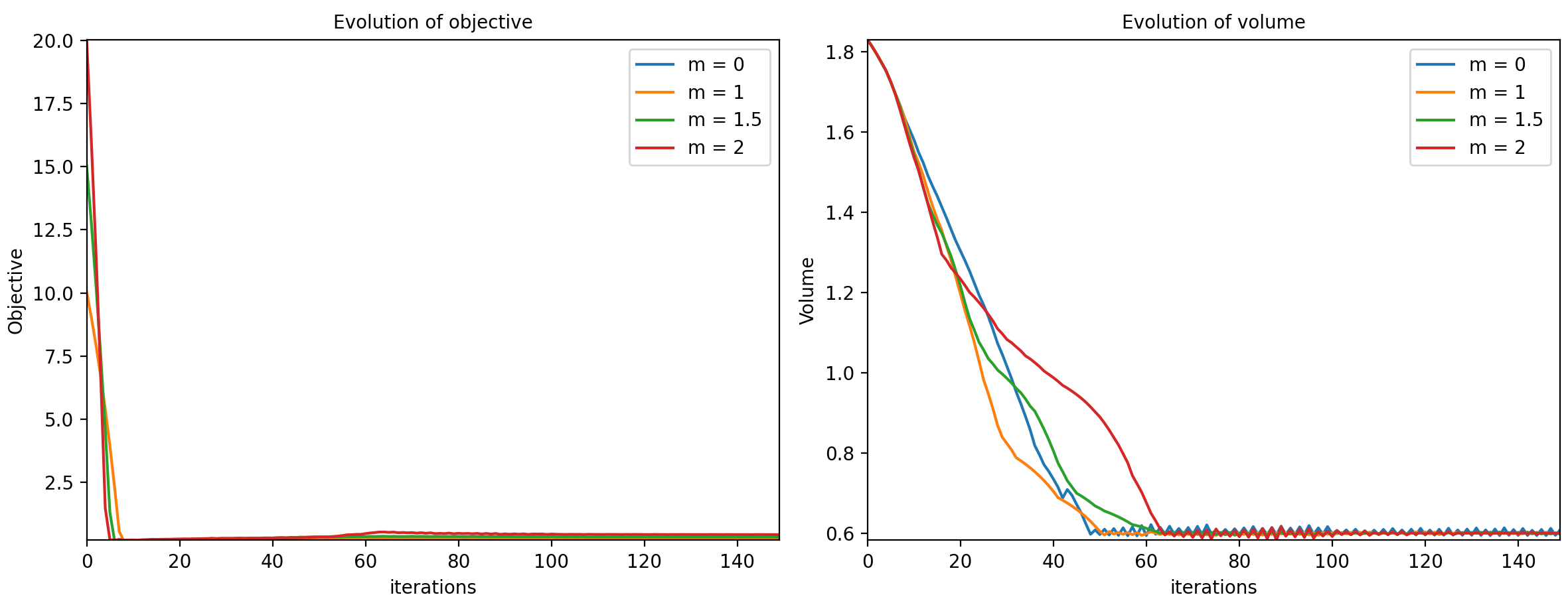}

\caption{\it Convergence histories for the experiments conducted with the cantilever shape optimization example of \cref{sec.canti}.}
\label{fig.cantihisto}
\end{figure}

We solve several instances of the distributionally robust problem \cref{eq.droaug} for the volume target $V_T = 0.6$ and various values of $m$ with the help of the constrained optimization algorithm introduced in \cite{feppon2020null}; 
the results are reported on \cref{fig.canti}, see \cref{fig.cantihisto} for the histories of the computations. 
Again, the nominal performance of the designs tends to deteriorate when the size of the parameter $m$ increases, see \cref{tab.cantitab}.

\begin{table}[!h]
\centering
	\begin{tabular}{|l|c|c|c|c|}
	\hline
	 Value of $m$  & $0$    &   $1$   & $1.5$   & $2$   \\
	\hline
	Nominal compliance & $0.0646956$    & $0.0746958$      & $0.08138$      &   $0.0963393$          \\
	\hline
	\end{tabular}
	\caption{\it Values of the nominal cost for the optimized cantilever beams of \cref{sec.canti}.}
\label{tab.cantitab}
\end{table}

\CD{
\begin{remark}
In both experiments, the resulting designs with $m=0$ do not correspond exactly with those 
that we would obtain by minimizing the compliance in the situation where only the nominal load is applied (i.e. by solving \cref{eq.nomoptpb}).
For example, in the instance of the cantilever test case of \cref{sec.canti} with $m=0$, we would intuitively expect to obtain a straight horizontal bar joining the regions $\Gamma_N$ and $\Gamma_D$.
Omitting the fact that, rigorously speaking, the duality result of \cref{prop.supWass} may not hold when $m>0$, this apparent inconsistency is due to the fact that, for fixed $\e >0$, the regularized Wasserstein distance $W_\e(\P,\Q)$ in \cref{eq.We} is only an approximation of the true Wasserstein distance $W(\P,\Q)$ in \cref{eq.Wass}. In particular, it is well-known in the literature (see e.g. \cite{feydy2020analyse}) that $W_\e(\P,\Q)$ is not a distance, so that the supremum in the left-hand side of \cref{eq.dualform} may not be attained when the probability law $\Q$ is the nominal law $\P$.  
\end{remark}
}

\CD{
\section{Conclusions and perspectives}\label{sec.concl}

\noindent To the best of our knowledge, this note is the first contribution devoted to distributionally robust formulations of shape and topology optimization problems featuring uncertain data. 
As we have seen, these consist in optimizing a design $h$ with respect to the worst mean value of a cost function $\calC(h,\xi)$ when the uncertain data $\xi$ follow a probability law $\Q$ which is ``close'' (in the sense of the entropy-regularized Wasserstein distance) to an observed nominal law $\P$. These very natural problems can be given convenient and tractable reformulations by leveraging
recent results from convex duality theory.
This work opens the way to multiple extensions, the perhaps most crucial of which being its extension to more realistic physical contexts where, for instance, the nominal law $\P$ is constructed from more than one observed sample, which entails a substantial increase in computational burden. 
It is also of utmost importance to compare the results of distributionally robust formulations with those of more classical robust formulations, like those mentioned in the introduction.
These questions, among others, will be addressed in a future and more complete work. 
}
\par\bigskip
\noindent \textbf{Conflict of interest.} On behalf of all authors, the corresponding author states that there is no conflict of interest. 

\par\medskip
\noindent \textbf{Replication of results.} The numerical example of \cref{sec.br} is tackled with the open-source finite element environment \texttt{FreeFem} \cite{hecht2012new}, 
and the precise source code used for the resolution is available on demand. 
The treatment of the numerical example of \cref{sec.canti} relies on minor adaptations to the open-source, educational implementation supplied with the article \cite{dapogny2022tuto}.  \par\medskip

\noindent \textbf{Acknowledgements.} The work of C.D. was partially supported by the project ANR-18-CE40-0013 SHAPO financed by the French Agence Nationale de la Recherche (ANR).

\CD{
\appendix
\section{Formal sketch of the proof of \cref{prop.supWass}}\label{app.proof}

\noindent 
At first, following a common practice in optimization theory, let us write the constrained supremum $A:= \sup_{W_\e(\P,\Q) \leq m} \int_\Xi f(\zeta) \:\d \Q(\zeta)$ as the unconstrained supremum 
of a quantity involving a Lagrange multiplier $\lambda$:
$$
A = \sup \limits_{\Q \in \calP(\Xi)} \inf\limits_{\lambda \geq 0 }\left\{ \int_\Xi f(\zeta) \:\d \Q(\zeta) + \lambda \left( m - W_\e(\P,\Q)\right) \right\}.
$$
Exchanging the above supremum and infimum, we obtain:
$$
A = \inf\limits_{\lambda \geq 0 } \sup \limits_{\Q \in \calP(\Xi)}  \left\{ \int_\Xi f(\zeta) \:\d \Q(\zeta) + \lambda \left( m - W_\e(\P,\Q)\right) \right\};
$$
the rigorous justification of this operation relies on the assumptions that $m>0$ and that $\sigma$ is small enough, 
which in brief guarantee that the maximization problem in the definition of $A$ is strictly feasible, 
see \cite{azizian2022regularization} for the details.

Now, the insertion of the definition \cref{eq.We} of $W_\e(\P,\Q)$ into this formula yields:
$$  A =  \inf\limits_{\lambda \geq 0 }\left\{ \lambda m + \sup \limits_{\Q \in \calP(\Xi)} \sup\limits_{\pi \in \calP(\Xi\times\Xi), \atop \pi_1 = \P, \pi_2 = \Q} \left\{ \int_\Xi f(\zeta) \:\d \Q(\zeta) - \lambda \int_{\Xi \times \Xi} c(\xi,\zeta) \:\d \pi(\xi,\zeta) - \lambda\e H(\pi) \right\} \right\}.$$
This formula features a maximization over all probability measures $\Q \in \calP(\Xi)$ and couplings $\pi \in \calP(\Xi \times \Xi)$ whose first and second marginals coincide with $\P$ and $\Q$ respectively; 
this is equivalent to a maximization over couplings $\pi \in \calP(\Xi \times \Xi)$ with $\P$ as first marginal and arbitrary second marginal. Hence, $A$ rewrites:
$$ A =  \inf\limits_{\lambda \geq 0 }\left\{ \lambda m + B(\lambda) \right\}, \text{ where } B(\lambda) := \sup\limits_{\pi \in \calP(\Xi\times\Xi) \atop \pi_1 = \P} \left\{ \int_{\Xi\times \Xi} \left( f(\zeta) - \lambda c(\xi,\zeta) \right) \:\d \pi(\xi,\zeta)  - \lambda \e H(\pi)\right\}.$$
The quantity $H(\pi)$ takes finite values only when $\pi$ is absolutely continuous with respect to $\pi_0$, i.e. when there exists a non negative integrable function $\alpha \in L^1(\Xi \times \Xi, \pi_0)$ 
such that $\pi = \alpha \pi_0$. Besides, from the definition of the reference coupling $\pi_0$, $\pi = \alpha \pi_0$ satisfies the marginal condition $\pi_1 = \P$ if and only if 
$$ \int_\Xi \alpha(\xi,\zeta) \d \nu_\xi(\zeta) = 1 \text{ for } \P\text{-a.e. } \xi \in \Xi.$$
Hence, $B(\lambda)$ rewrites:
\begin{equation}\label{eq.Blambda}
 B(\lambda) =  \sup\limits_{\alpha \in L^1(\Xi\times\Xi,\d\pi_0), \atop { \alpha \geq 0 , \atop \int_\Xi \alpha(\xi,\zeta) \d \nu_\xi(\zeta) = 1}} \left\{  \int_{\Xi\times \Xi}  \left( f(\zeta) - \lambda c(\xi,\zeta) - \lambda\e \log\alpha(\xi,\zeta) \right) \alpha(\xi,\zeta)\:\d \pi_0(\xi,\zeta) \right\}.
 \end{equation}
 
We now proceed to calculate $B(\lambda)$ explicitly, and
to this end, we rely on the first-order necessary condition for optimality in the above maximization problem.
According to the latter, there exists a function $\mu \in L^\infty(\Xi,\P)$ such that: 
$$ f(\zeta) - \lambda c(\xi,\zeta) - \lambda\e \log \alpha(\xi,\zeta) - \lambda\e + \mu(\xi) = 0 \text{ for } \pi_0\text{-almost every } (\xi,\zeta) \in \Xi \times \Xi.$$
Hence, we obtain the following expression for $\alpha(\xi,\zeta)$: 
\begin{equation*}\label{eq.alpha} 
\alpha(\xi,\zeta) = e^{\frac{f(\zeta) - \lambda c(\xi,\zeta) - \lambda \e + \mu(\xi)}{\lambda\e}},
\end{equation*}
and the expression for $\mu(\xi)$ is obtained from the constraint $\int_\Xi a(\xi,\zeta) \:\d \nu_\xi(\zeta) = 1$ for $\P$-a.e. $\xi \in \Xi$:
\begin{equation*}\label{eq.mu}
 \mu(\xi) = \lambda \e - \lambda \e \log\left(\int_\Xi e^{\frac{f(\zeta) - \lambda c(\xi,\zeta)}{\lambda\e}} \: \d\nu_\xi(\zeta) \right),
 \end{equation*}
so that eventually:
 $$\alpha(\xi,\zeta) = \left( \int_\Xi e^{\frac{f(\zeta) - \lambda c(\xi,\zeta)}{\lambda\e}} \:\d \nu_\xi(\zeta)  \right)^{-1} e^{\frac{f(\zeta) - \lambda c(\xi,\zeta)}{\lambda\e}}.$$ 
Inserting this last expression into \cref{eq.Blambda}, we arrive at:
 $$ B(\lambda) = \lambda \e \int_\Xi  \log \left(  \int_\Xi e^{\frac{f(\zeta) - \lambda c(\xi,\zeta)}{\lambda\e}} \:\d \nu_\xi(\zeta)  \right) \:\d \P(\xi),$$
whence the desired expression for $A$ is easily inferred.
 }
\bibliographystyle{siam}
\bibliography{./genbib.bib}

\end{document}